\documentclass[conference]{IEEEtran}
\makeatletter
\def\ps@headings{%
\def\@oddhead{\mbox{}\scriptsize\rightmark \hfil \thepage}%
\def\@evenhead{\scriptsize\thepage \hfil \leftmark\mbox{}}%
\def\@oddfoot{}%
\def\@evenfoot{}}
\makeatother 

\usepackage{cite}      % Written by Donald Arseneau
\usepackage{subfigure} % Written by Steven Douglas Cochran
\usepackage{url}       % Written by Donald Arseneau
\usepackage{algorithm}
\usepackage{array}
\usepackage{amsmath,amsfonts,amssymb,graphicx}
\usepackage{bm}
\usepackage{mathrsfs}
\usepackage{multicol}
\usepackage[end]{algpseudocode}
\usepackage{multirow}
\usepackage{dsfont}

\newtheorem{theorem}{\textbf{Theorem}}
\newtheorem{lemma}{\textbf{Lemma}}

\newtheorem{remark}{\textbf{Remark}}

\newcommand {\xA} {\mathcal{A}}

\newcommand{\bX} {\hat{\overline{\theta}}}

\newcommand{\qed}{\begin{flushright} $\blacksquare$ \end{flushright}}

\newcommand{\commnt}[1] {$//$ \textsc{#1} }
\newcommand{\xT} {\widetilde{T}_{i}}
\newcommand{\xI} {\widetilde{I}_{i}}
\newcommand{\xF} {\mathcal{F}}
\newcommand {\xa} {\mathbf{a}}
\newcommand{\xX} {\mathbf{X}}
\begin{document}
\title{Combinatorial Network Optimization with Unknown Variables: Multi-Armed Bandits with Linear Rewards}

\author{\IEEEauthorblockN{Yi Gai, Bhaskar Krishnamachari and Rahul Jain}
\IEEEauthorblockA{Ming Hsieh Department of Electrical Engineering\\
University of Southern California\\
Los Angeles, CA 90089, USA\\
Email: $\{$ygai, bkrishna, rahul.jain$\}$@usc.edu}}

\maketitle

\begin{abstract}

In the classic multi-armed bandits problem, the goal is to have a
policy for dynamically operating arms that each yield stochastic
rewards with unknown means. The key metric of interest is regret,
defined as the gap between the expected total reward accumulated by
an omniscient player that knows the reward means for each arm, and
the expected total reward accumulated by the given policy. The
policies presented in prior work have storage, computation and
regret all growing linearly with the number of arms, which is not
scalable when the number of arms is large. We consider in this work
a broad class of multi-armed bandits with dependent arms that yield
rewards as a linear combination of a set of unknown parameters. For
this general framework, we present efficient policies that are shown
to achieve regret that grows logarithmically with time, and
polynomially in the number of unknown parameters (even though the
number of dependent arms may grow exponentially). Furthermore, these
policies only require storage that grows linearly in the number of
unknown parameters. We show that this generalization is broadly
applicable and useful for many interesting tasks in networks that
can be formulated as tractable combinatorial optimization problems
with linear objective functions, such as maximum weight matching,
shortest path, and minimum spanning tree computations.

\end{abstract}

\section{Introduction}\label{sec:intro}

The problem of multi-armed bandits (MAB) is a classic one in
learning theory. In its simplest form, there are $N$ arms, each
providing stochastic rewards that are independent and identically
distributed over time, with unknown means. A policy is desired to
pick one arm at each time sequentially, to maximize the reward. MAB
problems capture a fundamental tradeoff between exploration and
exploitation; on the one hand, various arms should be explored in
order to learn their parameters, and on the other hand, the prior
observations should be exploited to gain the best possible immediate
rewards. MABs have been applied in a wide range of domains including
Internet advertising~\cite{Pandey,Rusmevichientong:2006} and
cognitive radio networks~\cite{Anandkumar:2010,Liu:Zhao}.

As they are fundamentally about combinatorial optimization in
unknown environments, one would indeed expect to find even broader
use of multi-armed bandits. However, we argue that a barrier to
their wider application in practice has been the limitation of the
basic formulation and corresponding policies, which generally treat
each arm as an independent entity. They are inadequate to deal with
many combinatorial problems of practical interest in which there are
large (exponential) numbers of arms. In such settings, it is
important to consider and exploit any structure in terms of
dependencies between the arms. We show in this work that when the
dependencies take a linear form, they can be handled tractably with
policies that have provably good performance in terms of regret as
well as storage and computation.

In this work, we formulate and consider the following general
multi-armed bandit problem. There is a vector $\xX$ of $N$ random
variables with unknown mean that are each instantiated in an i.i.d.
fashion over time. There is a finite (possibly exponentially large)
set of vector actions $\xa \in \xF$ from which any action can be
selected at each time. When action $\xa$ is performed, all elements
of $\xX$ that correspond to non-zero elements of $\xa$ are observed,
and a linear reward $\xa^T \xX$ is obtained. This generalization
captures a very broad class of combinatorial optimization problems
with linear objectives and unknown random coefficients.

A naive application of existing approaches for multi-armed bandits,
such as the well-known UCB1 index policy of Auer \emph{et
al.}~\cite{Auer:2002}, for this problem would yield poor performance
scaling in terms of regret, storage, and computation. This is
because these approaches are focused on maintaining and computing
quantities based on arm-specific observations and do not exploit
potential dependencies between them. In this work, we instead
propose smarter policies that explicitly take into account the
linear form of the dependencies and base all storage and
computations on the unknown variables directly, rather than the
arms. As we shall show, this saves not only on storage and
computation, but also substantially reduces the regret.

Specifically, we first present a novel single-arm selection policy
for Learning with Linear Rewards (LLR) requires only $O(N)$ storage,
and yields a regret that grows essentially~\footnote{This is a
simplification of our key result in section~\ref{sec:regret} which
gives a tighter expression for the bound on regret that applies
uniformly over time, not just asymptotically.} as $O(N^4 \ln n)$,
where $n$ is the time index. We also discuss how this policy can be
modified in a straightforward manner while maintaining the same
performance guarantees when the problem is one of cost minimization
rather than reward maximization. A key step in these policies we
propose is the solving of a deterministic combinatorial optimization
with a linear objective. While this is NP-hard in general (as it
includes 0-1 integer linear programming), there are still many
special-case combinatorial problems of practical interest which can
be solved in polynomial time. For such problems, the policy we
propose would thus inherit the property of polynomial computation at
each step.

We also present in this paper a more general K-arm formulation, in
which the policy is allowed to pick $K \geq 1$ different actions at
each time. We show how the single-arm policy can be readily extended
to handle this and present the regret analysis for this case as
well.

Through several concrete examples, we show the applicability of our
general formulation of multi-armed bandits with linear rewards to
combinatorial network optimization. These include maximum weight
matching in bipartite graphs (which is useful for user-channel
allocations in cognitive radio networks), as well as shortest path,
and minimum spanning tree computation. The examples we present are
far from exhausting the possible applications of the formulation and
the policies we present in this work --- there are many other
linear-objective network optimization problems~\cite{Ahuja:1993,
Vygen:2008}. Our framework, for the first time, allows these
problems to be solved in stochastic settings with unknown random
coefficients, with provably efficient performance.

We expect that our work will also find practical application in
other fields where such linear combinatorial optimization problems
arise naturally, such as algorithmic economics, data mining,
finance, operations research and industrial engineering.

This paper is organized as follows. We first provide a survery of
related work in section \ref{sec:related}. We then give a formal
description of the multi-armed bandits with linear rewards problem
we solve in section \ref{sec:formulation}. In section
\ref{sec:general}, we present our LLR policy and show that it
requires only polynomial storage and polynomial computation per time
period. We present the novel analysis of the regret of this policy
in section \ref{sec:regret} and point out how this analysis
generalizes known results on MAB. In section \ref{sec:example}, we
discuss examples and applications of maximum weight matching,
shortest path, and minimum spanning tree computations to show that
our policy is widely useful for various interesting applications in
networks with the tractable combinatorial optimization formulation
with linear objective functions. Section~\ref{sec:simulation} shows
the numerical simulation results. We show an extension of our policy
for choosing $K$ largest values in section \ref{sec:klargest}.
Finally, we conclude with a summary of our contribution and point
out avenues for future work in section \ref{sec:conclusion}.

\section{Related Work}\label{sec:related}
%need to expand this section greatly to be more comprehensive

Lai and Robbins~\cite{Lai:Robbins} wrote one of the earliest papers
on the classic non-Bayesian infinite horizon multi-armed bandit
problem. Assuming K independent arms, each generating rewards that
are i.i.d. over time from a given family of distributions with an
unknown real-valued parameter, they presented a general policy that
provides expected regret that is $O(K \log n)$, i.e. linear in the
number of arms and asymptotically logarithmic in n. They also show
that this policy is order optimal in that no policy can do better
than $\Omega(K \log n)$. Anantharam \emph{et al.}~\cite{Anantharam}
extend this work to the case when $M$ simultaneous plays are
allowed. The work by Agrawal~\cite{Agrawal:1995} presents easier to
compute policies based on the sample mean that also has
asymptotically logarithmic regret. However, their policies need not
be directly applied to our problem formulation in this paper, which
involves combinatorial arms that cannot be characterized by a single
parameter.

Our work is influenced by the paper of Auer \emph{et
al.}~\cite{Auer:2002} that considers arms with non-negative rewards
that are i.i.d. over time with an arbitrary un-parameterized
distribution that has the only restriction that it have a finite
support. Further they provide a simple policy (referred to as UCB1),
which achieves logarithmic regret uniformly over time, rather than
only asymptotically. However, their work does not exploit potential
dependencies between the arms. As we show in this paper, a direct
application of their UCB1 policy therefore performs poorly for our
problem formulation.

There are also some recent works to propose decentralized policies
for the multi-armed bandit problem. Liu and Zhao~\cite{Liu:Zhao},
and Anandkumar \emph{et al.}~\cite{Anandkumar:2010} have both
developed policies for the problem of $M$ distributed players
operating $N$ independent arms.

%There are also some recent work to propose the decentralized
%policies for multi-armed bandit problem. Liu and
%Zhao~\cite{Liu:Zhao} formulated the problem of $M$ distributed
%players competing for $N$ independent arms, (each generating rewards
%that are i.i.d. over time from a given family of distributions with
%an unknown real-valued parameter). They present a policy that
%achieves asymptotically logarithmic regret with respect to time.
%Anandkumar \emph{et al.}~\cite{Anandkumar:2010} solves the similar
%decentralized problem based on the single-user policy proposed in
%\cite{Auer:2002}. They propose two decentralized learning and access
%policies, one assuming the total number of users is know and the
%other removes this requirement. They prove that the first policy
%achieves logarithmic regret and the second achieves nearly
%logarithmic regret.

While these above key papers and many others have focused on
independent arms, there have been some works treating dependencies
between arms. The paper by Pandey \emph{et al.}~\cite{Pandey}
divides arms into clusters of dependent arms (in our case there
would be only one such cluster consisting of all the arms). Their
model assumes that each arm provide only binary rewards, and in any
case, they do not present any theoretical analysis on the expected
regret. Ortner~\cite{Ortner:2010} proposes to use an additional arm
color, to utilize the given similarity information of different arms
to improve the upper bound of the regret. They assume that the
difference of the mean rewards of any two arms with the same color
is less than a predefined parameter $\delta$, which is known to the
user. This is different from the linear reward model in our paper.

Mersereau \emph{et al.}~\cite{Mersereau:2008} consider a bandit
problem where the expected reward is defined as a linear function of
an random variable, and the prior distribution is known. They show
the upper bound of the regret is $O(\sqrt{n})$ and the lower bound
of the regret is $\Omega(\sqrt{n})$. Rusmevichientong and Tsitsiklis
\cite{Rusmevichientong:2010} extend~\cite{Mersereau:2008} to the
setting where the reward from each arm is modeled as the sum of a
linear combination of a set of unknown static random numbers and a
zero-mean random variable that is i.i.d. over time and independent
across arms. The upper bound of the regret is shown to be $O(N
\sqrt{n})$ on the unit sphere and $O(N \sqrt{n} \log^{3/2} n)$ for a
compact set, and the lower bound of regret is $\Omega(N \sqrt{n})$
for both cases. The linear models in these works are different from
our paper in which the reward is expressed as a linear combination
as a set of random processes. Also,~\cite{Mersereau:2008}
and~\cite{Rusmevichientong:2010} assume that only the reward is
observed at each time. In our work, we assume that the random
variables corresponding to non-zero action components are observed
at each time (from which the reward can be inferred).

Both~\cite{Auer:confidence} and~\cite{Dani:2008} consider linear
reward models that are more general than ours, but also under the
assumption that only the reward is observed at each time. Auer
\cite{Auer:confidence} presents a randomized policy which requires
storage and computation to grow linearly in the number of arms. This
algorithm is shown to achieve a regret upper bound of
$O(\sqrt{N}\sqrt{n} \log^{\frac{3}{2}}(n|\xF|))$. Dani \emph{et al.}
\cite{Dani:2008} develop another randomized policy for the case of a
compact set of arms, and show the regret is upper bounded by $O(N
\sqrt{n} \log^{3/2} n)$ for sufficiently large $n$ with high
probability, and lower bounded by $\Omega(N \sqrt{n})$. They also
show that when the difference in costs (denoted as $\Delta$) between
the optimal and next to optimal decision among the extremal points
is greater than zero, the regret is upper bounded by
$O(\frac{N^2}{\Delta} \log^{3} n)$ for sufficiently large $n$ with
high probability. To our best knowledge, ours is the first paper to
consider linear rewards with observation of the random variables
corresponding to non-zero action components. We present a
deterministic policy with a deterministic combinatorial linear
optimization problem finite time bound of regret which grows $O(N^4
\log n)$, i.e., polynomially in the number of unknown random
variables and strictly logarithmically in time.

%They also show that when the arms are from a finite set, the regret
%is upper bounded by $O(\frac{N^2}{\delta} \log^{3} n)$ for
%sufficiently large $n$ with high probability.

Our work in this paper is an extension of our recent work which
introduced combinatorial multi-armed bandits \cite{Gai:2010}. The
formulation in \cite{Gai:2010} has the restriction that the reward
is generated from a matching in a bipartite graph of users and
channels. Our work in this paper generalizes this to a broader
formulation with linear reward, where the action vector is from a
finite set.

\section{Problem Formulation}\label{sec:formulation}

Now we define the problem of multi-armed bandits with linear rewards
that we solve in this paper. We consider a discrete time system with
$N$ unknown random processes $X_i(n), 1 \leq i \leq N $, where time
is indexed by $n$. We assume that $X_i(n)$ evolves as an i.i.d.
random process over time, with the only restriction that its
distribution have a finite support. Without loss of generality, we
normalize $X_i(n) \in [0,1]$. We do not require that $X_i(n)$ be
independent across $i$. This random process is assumed to have a
mean $\theta_{i} = E[X_i]$ that is unknown to the users. We denote
the set of all these means as $\Theta = \{\theta_{i}\}$.

At each decision period $n$ (also referred to interchangeably as
time slot),  an $N$-dimensional action vector $\xa(n)$, representing
an arm, is selected under a policy $\pi(n)$ from a finite set $\xF$.
We assume $a_i(n) \geq 0$ for all $1 \leq i \leq N$. When a
particular $\xa(n)$ is selected, only for those $i$ with $a_i(n)
\neq 0$, the value of $X_i(n)$ is observed . We denote $\xA_{\xa(n)}
= \{i: a_i(n) \neq 0, 1 \leq i \leq N \}$, the index set of all
$a_i(n) \neq 0$ for an arm $\xa$. We treat each $\xa(n) \in \xF$ as
an arm. The reward is defined as:
\begin{equation}
  R_{\xa(n)}(n) = \sum\limits_{i = 1}^N a_i(n) X_i(n).
\end{equation}

When a particular action/arm $a(n)$ is selected, the random
variables corresponding to non-zero components of a(n) are
revealed\footnote{As noted in the related work, this is a key
assumption in our work that differentiates it from other prior work
on linear dependent-arm bandits \cite{Auer:confidence},
\cite{Dani:2008}. This is a very reasonable assumption in many
cases, for instance, in the combinatorial network optimization
applications we discuss in section \ref{sec:example}, it corresponds
to revealing weights on the set of edges selected at each time.},
i.e., the value of $X_i(n)$ is observed for all i such that $\xa(n)
\neq 0$.

%We define the policy $\pi(n)$ at each time to be a map from the
%observation history $\{O_k\}_{k=1}^{n-1}$ to a vector of channels
%$o(n)$ to be selected at period $n$, where user $i$ selects channel
%$o_i(n)$. Then the observation history $\{O_k\}_{k=1}^{n-1}$ in turn
%can be expressed as $\{o_i(k), X_{i, o_i(k)}(k)\}_{1 \leq i \leq M,
%1 \leq k < n}$.

%%%%%%%%should be in related work%%%%%%%%%%%%%%
%This kind of problem, where the desired goal is to develop a
%sequential policy to make a selection among multiple choices, each
%offering stochastic rewards derived from a distribution with an
%unknown parameter, is traditionally formulated as an infinite
%horizon non-Bayesian multi-armed bandit (see~\cite{Lai:Robbins,
%Anantharam, Agrawal:1995, Auer:2002}).

We evaluate policies with respect to \emph{regret}, which is defined
as the difference between the expected reward that could be obtained
by a genie that can pick an optimal arm at each time, and that
obtained by the given policy. Note that minimizing the regret is
equivalent to maximizing the rewards. Regret can be expressed as:
\begin{equation}
 \mathfrak{R}^\pi_n (\Theta) = n \theta^*  - E^\pi[ \sum \limits_{t = 1}^n R_{\pi(t)}(t) ],
\end{equation}
where $\theta^* = \max\limits_{\xa \in\xF} \sum\limits_{i = 1}^N a_i
\theta_i$, the expected reward of an optimal arm. For the rest of
the paper, we use $*$ as the index indicating that a parameter is
for an optimal arm. If there is more than one optimal arm exist, $*$
refers to any one of them.

Intuitively, we would like the regret $\mathfrak{R}^\pi_n (\Theta)$
to be as small as possible. If it is sub-linear with respect to time
$n$, the time-averaged regret will tend to zero and the maximum
possible time-averaged reward can be achieved. Note that the number
of arms $|\xF|$ can be exponential in the number of unknown random
variables $N$.

%Different from most prior work, although the number of unknown
%variables are $N$, the number of arms, denoted as $|\xF|$, may not
%be polynomial in $N$, and could be exponential in $N$.
%For example,
%in the maximum bipartite matching example in section
%\ref{sec:example}, the number of arms grows as the permutation
%$P(N,M)$. Also, note that the reward obtained from each arm is
%dependent across arms that share common components.

%There are many practical problems fit the model of this problem,
%such as the multiuser channel allocation problem in cognitive
%networks and shortest path problem to minimum the delay. We will
%present in Section \ref{sec:example} in details.

\section{Policy Design}\label{sec:general}

\subsection{A Naive Approach}\label{sec:naive}

A straightforward, relatively naive approach to solving the
multi-armed bandits with linear regret problem that we defined is to
use the UCB1 policy given by Auer~\emph{et al.}~\cite{Auer:2002}.
For UCB1, the arm that maximizes $\hat{Y}_k + \sqrt{ \frac{2 \ln n
}{ m_k } }  $ will be selected at each time slot, where $\hat{Y}_k$
is the mean observed reward on arm $k$, and $m_k$ is the number of
times that arm $k$ has been played. This approach essentially
ignores the dependencies across the different arms, storing observed
information about each arm independently, and making decisions based
on this information alone.

Auer~\emph{et al.}~\cite{Auer:2002} showed the following policy
performance for regret upper bound as:

\theorem\label{cl:classic} The expected regret under UCB1 policy is
at most

\begin{equation}
 \left[ 8\sum\limits_{k: \theta_k < \theta^*} ( \frac{\ln n}{\Delta_k} ) \right] + (1+ \frac{\pi^2}{3})
 ( \sum\limits_{k: \theta_k < \theta^* } \Delta_k)
\end{equation}
where $\Delta_k  = \theta^* - \theta_k$, $\theta_k = \sum \limits_{i
\in \xA_k} a_i \theta_{i}$.

\textit{Proof:} See~\cite[Theorem 1]{Auer:2002}. \qed

Note that UCB1 requires storage that is linear in the number of arms
and yields regret growing linearly with the number of arms. In a
case where the number of arms grow exponentially with the number of
unknown variables, both of these are highly unsatisfactory.

Intuitively, UCB1 algorithm performs poorly on this problem because
it ignores the underlying dependencies. This motivates us to propose
a sophisticated policy which more efficiently stores observations
from correlated arms and exploits the correlations to make better
decisions.

\subsection{A new policy}\label{sec:generalAlg}
%%%%%%%%%%%%%%%%%%%%%%%%%%%%%%%%%%%%%%%%%%%%%%%%%%%%%%%%%%%%%%%%%%%%%%%%%%%%%%%%%%%%%%%%%%%%%%%%%%%%%
%%%%%%%%%%%%%%%%%%%%%%%%%%%%%%%%%%%%%%%%%%%%%%%%%%%%%%%%%%%%%%%%%%%%%%%%%%%%%%%%%%%%%%%%%%%%%%%%%%%%%
Our proposed policy, which we refer to as ``learning with linear
rewards" (LLR), is shown in Algorithm \ref{alg:general}.

%%%%%%%%%%%%%%%%%%%%%%
%%%%%%%%%%%%%%%%%%%%%%%
\begin{algorithm} [ht]
\caption{Learning with Linear Rewards (LLR)} \label{alg:general}

\begin{algorithmic}[1]
%----------------------------------------------------------------
\State \commnt{ Initialization}

\State If $\max\limits_\xa |\xA_\xa|$ is known, let $L =
\max\limits_\xa |\xA_\xa|$; else, $L = N$;

\For {$p = 1$ to $N$}
    \State $n = p$;
    \State Play any arm $\xa$ such that $p \in \xA_\xa$;
    \State Update $(\hat{\theta}_{i})_{1 \times N}$, $(m_{i})_{1 \times N}$
    accordingly;
\EndFor

\State \commnt{Main loop}

\While {1}
    \State $n = n + 1$;
    \State Play an arm $\xa$ which solves the maximization problem
    \begin{equation}
    \label{equ:choosemax}
    \xa = \arg\max\limits_{\xa \in \xF} \sum\limits_{i \in \xA_\xa} a_i \left( \hat{ \theta }_{i} + \sqrt{ \frac{ (L+1) \ln n }{ m_{i}}
    } \right);
    \end{equation}
    \State Update $(\hat{\theta}_{i})_{1 \times N}$, $(m_{i})_{1 \times N}$
    accordingly;
\EndWhile
\end{algorithmic}
\end{algorithm}
%%%%%%%%%%%%%%%%%%%%%%%
%%%%%%%%%%%%%%%%%%%%%%%
\begin{table}[htbp]
\centering \normalsize
\begin{tabular}{|l l|}
\hline
&\\
 $N$ : & number of random variables.\\
 $\xa$ : & vectors of coefficients, defined on set $\xF$; \\
  & we map each $\xa$ as an arm.\\
 $\xA_\xa$: & $\{i: a_i \neq 0, 1 \leq i \leq N \}$.\\
 $*$ : & index indicating that a parameter is for an \\
 & optimal arm. \\
 $m_{i}$: & number of times that $X_i$ has been observed \\
  & up to the current time slot.\\
 $\hat{ \theta }_{i}$: & average (sample mean) of all the observed \\
  & values of $X_i$ up to the current time slot.\\
 & Note that $\mathbb{E}[\hat{ \theta }_{i}(n)] = \theta_{i}$.\\
 $\bX_{i, m_i}$: & average (sample mean) of all the observed \\
 & values of $X_i$ when it is observed $m_{i}$ times.\\
 $\Delta_\xa$: & $R^* - R_{\xa}$ .\\
 $\Delta_{\min}$: & $\min\limits_{\xa \neq \xa^*} \Delta_\xa$.\\
 $\Delta_{\max}$: & $\max\limits_{\xa \neq \xa^*} \Delta_\xa$.\\
 $T_\xa (n)$: & number of times arm $\xa$ has been played\\
 &in the first $n$ time slots.\\
 $a_{\max}$: & $\max\limits_{\xa \in \xF} \max\limits_{i} a_i$.\\
\hline
\end{tabular}
\caption{Notation} \label{tab:notation}
\end{table}

Table \ref{tab:notation} summarizes some notation we use in the
description and analysis of our algorithm.

The key idea behind this algorithm is to store and use observations
for each random variable, rather than for each arm as a whole. Since
the same random variable can be observed while operating different
arms, this allows exploitation of information gained from the
operation of one arm to make decisions about a dependent arm.

We use two $1$ by $N$ vectors to store the information after we play
an arm at each time slot. One is $(\hat{ \theta }_{i})_{1 \times N}$
in which $\hat{ \theta }_{i}$ is the average (sample mean) of all
the observed values of $X_i$ up to the current time slot (obtained
through potentially different sets of arms over time). The other one
is $(m_{i})_{1 \times N}$ in which $m_{i}$ is the number of times
that $X_i$ has been observed up to the current time slot.

At each time slot $n$, after an arm $\xa(n)$ is played, we get the
observation of $X_{i}(n)$ for all $i \in \xA_{\xa(n)}$. Then $(\hat{
\theta }_{i})_{1 \times N}$ and $(m_{i})_{1 \times N}$ (both
initialized to 0 at time 0) are updated as follows:
\begin{equation}
 \hat{ \theta }_{i}(n) = \left\{ \begin{array}
  {l@{\quad,\quad}l}
  \frac{\hat{ \theta }_{i}(n-1) m_{i}(n-1) + X_{i}(n)}{ m_{i}(n-1) +1} & \text{if } i \in \xA_{\xa(n)}\\
  \hat{ \theta }_{i}(n-1) & \text{else} \\
  \end{array} \right.
\end{equation}
\begin{equation}
  m_{i}(n) = \left\{ \begin{array}
  {l@{\quad,\quad}l}
  m_{i}(n-1)+1 & \text{if } i \in \xA_{\xa(n)}\\
  m_{i}(n-1) & \text{else} \\
  \end{array} \right.
\end{equation}

Note that while we indicate the time index in the above updates for
notational clarity, it is not necessary to store the matrices from
previous time steps while running the algorithm.

LLR policy requires storage linear in $N$. In section
\ref{sec:regret}, we will present the analysis of the upper bound of
regret, and show that it is polynomial in $N$ and logarithmic in
time. Note that the maximization problem (\ref{equ:choosemax}) needs
to be solved as the part of LLR policy. It is a deterministic linear
optimal problem with a feasible set $\xF$ and the computation time
for an arbitrary $\xF$ may not be polynomial in $N$. As we show in
Section \ref{sec:example}, that there exists many practically useful
examples with polynomial computation time.

\section{Analysis of Regret}\label{sec:regret}

Traditionally, the regret of a policy for a multi-armed bandit
problem is upper-bounded by analyzing the expected number of times
that each non-optimal arm is played, and the summing this
expectation over all non-optimal arms. While such an approach will
work to analyze the LLR policy too, it turns out that the
upper-bound for regret consequently obtained is quite loose, being
linear in the number of arms, which may grow faster than
polynomials. Instead, we give here a tighter analysis of the LLR
policy that provides an upper bound which is instead polynomial in
$N$ and logarithmic in time. Like the regret analysis
in~\cite{Auer:2002}, this upper-bound is valid for finite $n$.

\theorem\label{c1:upperbound} The expected regret under the LLR
policy is at most
\begin{equation}
\left[\frac{ 4 a_{\max}^2 L^2 (L+1) N \ln n }{ \left( \Delta_{\min}
\right)^2 } + N + \frac{\pi^2}{3} L N \right] \Delta_{\max}.
\end{equation}

To proof Theorem \ref{c1:upperbound}, we use the inequalities as
stated in the Chernoff-Hoeffding bound \cite{Pollard}.

\lemma[Chernoff-Hoeffding bound \cite{Pollard}] \label{lemma:1}
$X_1, \ldots, X_n$ are random variables with range $[0,1]$, and
$E[X_t|X_1, \ldots, X_{t-1}] = \mu$, $\forall 1 \leq t \leq n$.
Denote $S_n = \sum X_i$. Then for all $a \geq 0$
\begin{equation}
 \begin{split}
  Pr\{S_n \geq n \mu + a\} & \leq e^{-2 a^2 /n}\\
  Pr\{S_n \leq n \mu - a\} & \leq e^{-2 a^2 /n}
 \end{split}
\end{equation}

\begin{IEEEproof}[Proof of Theorem \ref{c1:upperbound}] Denote $C_{t, m_i}$ as $\sqrt{ \frac{ (L+1) \ln t }{ m_i} }$.
We introduce $\xT(n)$ as a counter after the initialization period.
It is updated in the following way:

At each time slot after the initialization period, one of the two
cases must happen: (1) an optimal arm is played; (2) a non-optimal
arm is played. In the first case, $(\xT(n))_{1 \times N}$ won't be
updated. When an non-optimal arm $\xa(n)$ is picked at time $n$,
there must be at least one $i \in \xA_\xa$ such that $i = arg
\min\limits_{ j \in \xA_\xa } m_{j}$. If there is only one such arm,
$\xT(n)$ is increased by $1$. If there are multiple such arms, we
arbitrarily pick one, say $i'$, and increment $\widetilde{T}_{i'}$
by $1$.

Each time when a non-optimal arm is picked, exactly one element in
$(\xT(n))_{1 \times N}$ is incremented by $1$. This implies that the
total number that we have played the non-optimal arms is equal to
the summation of all counters in $(\xT(n))_{1 \times N}$. Therefore,
we have:
\begin{equation}
 \label{equ:f1}
 \sum\limits_{\xa: \xa \neq \xa^*} \mathbb{E}[T_\xa(n)] = \sum\limits_{i = 1}^N
 \mathbb{E}[\xT(n)].
\end{equation}

Also note for $\xT(n)$, the following inequality holds:
\begin{equation}
 \label{equ:f2}
 \xT(n) \leq m_{i}(n), \forall 1 \leq i \leq N.
\end{equation}

Denote by $\xI(n)$ the indicator function which is equal to $1$ if
$\xT(n)$ is added by one at time $n$. Let $l$ be an arbitrary
positive integer. Then:

\begin{equation} \label{equ:10}
\begin{split}
 \xT(n) & = \sum\limits_{t = N+1}^n \mathds{1} \{ \xI(t)=1\} \\
 & \leq l + \sum\limits_{t = N+1}^n \mathds{1} \{ \xI(t)=1 , \xT(t-1) \geq l \} \\
\end{split}
\end{equation}
where $\mathds{1}(x)$ is the indicator function defined to be 1 when
the predicate $x$ is true, and 0 when it is false. When $\xI(t) =
1$, a non-optimal arm $\xa(t)$ has been picked for which $m_{i} =
\min\limits_j \{m_j: \forall j \in \xA_{\xa(t)}\}$. We denote this
arm as $\xa(t)$ since at each time that $\xI(t) = 1$, we could get
different arms. Then,
\begin{equation} \label{equ:11}
\begin{split}
 \xT(n) & \leq l + \sum\limits_{t = N+1}^n \mathds{1}\{ \sum\limits_{j \in \xA_{\xa^*}} a_j^* (\bX_{j, m_j(t-1)} + C_{t-1, m_j(t-1)} )  \\
 \leq &  \sum\limits_{j \in \xA_{\xa(t)}} a_j(t) (\bX_{j, m_j(t-1)} + C_{t-1, m_j(t-1)} ), \xT(t-1) \geq l \} \\
 & \leq l + \sum\limits_{t = N}^n \mathds{1}\{ \sum\limits_{j \in \xA_{\xa^*}} a_j^* (\bX_{j, m_j(t)} + C_{t, m_j(t)} ) \\
 & \quad \leq \sum\limits_{j \in \xA_{\xa(t)}} a_j(t) (\bX_{j, m_j(t)} + C_{t, m_j(t)} ), \xT(t) \geq l
 \}.
\end{split}
\end{equation}
Note that $l \leq \xT(t)$ implies,
\begin{equation} \label{equ:13}
 l \leq \xT(t) \leq m_{j}(t), \forall j \in \xA_{\xa(t)}.
\end{equation}

%Then we have,
%\begin{equation}
%\begin{split}
% \xT(n) & \leq l + \sum\limits_{t = N}^n \{ \min\limits_{0 < m_i \leq t, \forall i \in \xA_{\xa^*}  }   \sum\limits_{i \in \xA_{\xa*}} a_i^* (\bX_{i, m_i} + C_{t, m_i} ) \nonumber\\
% & \leq \max\limits_{l \leq m_j \leq t, \forall j \in \xA_{\xa(t)} } \sum\limits_{j \in \xA_{\xa(t)}} a_j(t) (\bX_{i, m_j} + C_{t, m_j} ) \}\\
% & \leq l + \sum\limits_{t = 1}^{\infty}  \sum\limits_{m_i = 1, \forall i \in \xA_{\xa^*}}^t \sum\limits_{m_j = 1, \forall j \in \xA_{\xa(t)}}^t \\
% & \quad \{ \sum\limits_{i \in \xA_{\xa^*}} a_i^* (\bX_{i, m_i(t)} + C_{t, m_i(t)} ) \\
% & \quad \leq \sum\limits_{j \in \xA_{\xa(t)}} a_j(t) (\bX_{i, m_j} + C_{t, m_j} ) \}
%\end{split}
%\end{equation}

\begin{equation} \label{equ:12}
\begin{split}
 \xT(n) & \leq l + \sum\limits_{t = N}^n \mathds{1} \{ \min\limits_{0 < m_{h_1}, \ldots, m_{h_{|\xA_{\xa*}|}} \leq t }   \sum\limits_{j = 1}^{|\xA_{\xa*}|} a_{h_j}^* (\bX_{h_j, m_{h_j}} + C_{t, m_{h_j}} ) \\
 & \quad \leq \max\limits_{l \leq m_{p_1}, \ldots, m_{p_{|\xA_{\xa(t)}|}} \leq t} \sum\limits_{j = 1}^{|\xA_{\xa(t)}|} a_{p_j}(t) (\bX_{p_j, m_{p_j}} + C_{t, m_{p_j}} ) \}\\
 & \leq l + \sum\limits_{t = 1}^{\infty} \sum\limits_{m_{h_1} = 1}^{t} \dots \sum\limits_{m_{h_{|\xA^*|}} = 1}^{t}
 \sum\limits_{m_{p_1} = l}^{t} \dots \sum\limits_{m_{p_{|\xA_{\xa(t)}|}} = l}^{t}    \\
 & \quad\quad \mathds{1}\{\sum\limits_{j = 1}^{|\xA_{\xa*}|} a_{h_j}^* (\bX_{h_j, m_{h_j}} + C_{t, m_{h_j}} ) \\
 & \quad\quad \leq \sum\limits_{j = 1}^{|\xA_{\xa(t)}|} a_{p_j}(t) (\bX_{p_j, m_{p_j}} + C_{t, m_{p_j}} )  \}
\end{split}
\end{equation}
where $h_j$ ($1 \leq j \leq |\xA_{\xa*}|$) represents the $j$-th
element in $\xA_{\xa*}$ and $p_j$ ($1 \leq j \leq |\xA_{\xa(t)}|$)
represents the $j$-th element in $\xA_{\xa(t)}$.

$\sum\limits_{j = 1}^{|\xA_{\xa*}|} a_{h_j}^* (\bX_{h_j, m_{h_j}} +
C_{t, m_{h_j}}) \leq \sum\limits_{j = 1}^{|\xA_{\xa(t)}|} a_{p_j}(t)
(\bX_{p_j, m_{p_j}} + C_{t, m_{p_j}} )$ means that at least one of
the following must be true:
\begin{equation}
\label{equ:p1}
 \sum\limits_{j = 1}^{|\xA_{\xa*}|} a_{h_j}^* \bX_{h_j, m_{h_j}} \leq  R^* - \sum\limits_{j = 1}^{|\xA_{\xa*}|} a_{h_j}^* C_{t, m_{h_j}},
\end{equation}
\begin{equation}
\label{equ:p2}
 \sum\limits_{j = 1}^{|\xA_{\xa(t)}|} a_{p_j}(t) \bX_{p_j, m_{p_j}} \geq R_{\xa(t)} + \sum\limits_{j = 1}^{|\xA_{\xa(t)}|} a_{p_j}(t) C_{t, m_{p_j}},
\end{equation}
\begin{equation}
 R^* < R_{\xa(t)}  + 2 \sum\limits_{j = 1}^{|\xA_{\xa(t)}|} a_{p_j}(t) C_{t, m_{p_j}}.
\end{equation}
Now we find the upper bound for $Pr\{ \sum\limits_{j =
1}^{|\xA_{\xa*}|} a_{h_j}^* \bX_{h_j, m_{h_j}} \leq  R^* -
\sum\limits_{j = 1}^{|\xA_{\xa*}|} a_{h_j}^* C_{t, m_{h_j}} \}$.

We have:

\begin{equation}
\begin{split}
 & Pr\{ \sum\limits_{j = 1}^{|\xA_{\xa*}|} a_{h_j}^* \bX_{h_j, m_{h_j}} \leq  R^*
   - \sum\limits_{j = 1}^{|\xA_{\xa*}|} a_{h_j}^* C_{t, m_{h_j}} \} \nonumber\\
 & = Pr\{ \sum\limits_{j = 1}^{|\xA_{\xa*}|} a_{h_j}^* \bX_{h_j, m_{h_j}}
   \leq \sum\limits_{j = 1}^{|\xA_{\xa*}|} a_{h_j}^* \theta_{h_j}
     - \sum\limits_{j = 1}^{|\xA_{\xa*}|} a_{h_j}^* C_{t, m_{h_j}}\}  \\
 & \leq Pr\{ \text{At least one of the following must hold:} \\
 & \qquad \qquad a_{h_1}^* \bX_{h_1, m_{h_1}} \leq a_{h_1}^* \theta_{h_1} - a_{h_1}^* C_{t, m_{h_1}}, \\
 & \qquad \qquad a_{h_2}^* \bX_{h_2, m_{h_2}} \leq a_{h_2}^* \theta_{h_2} - a_{h_2}^* C_{t, m_{h_2}}, \\
 & \qquad \qquad \qquad \qquad \vdots \\
 & \qquad \qquad a_{h_{|\xA_{\xa*}|}}^* \bX_{h_1, m_{h_{|\xA_{\xa*}|}}} \leq a_{h_{|\xA_{\xa*}|}}^* \theta_{h_{|\xA_{\xa*}|}} \\
 & \qquad \qquad\qquad \qquad- a_{h_{|\xA_{\xa*}|}}^* C_{t, m_{h_{|\xA_{\xa*}|}}} \} \\
 & \leq \sum\limits_{j = 1}^{|\xA_{\xa*}|} Pr\{ a_{h_j}^* \bX_{h_j, m_{h_j}} \leq a_{h_j}^* \theta_{h_j} - a_{h_j}^* C_{t, m_{h_j}}
 \} \\
 & = \sum\limits_{j = 1}^{|\xA_{\xa*}|} Pr\{ \bX_{h_j, m_{h_j}} \leq \theta_{h_j} - C_{t, m_{h_j}}
 \}.
\end{split}
\end{equation}

$\forall 1 \leq j \leq {|\xA_{\xa*}|}$, applying the
Chernoff-Hoeffding bound stated in Lemma \ref{lemma:1}, we could
find the upper bound of each item in the above equation as,
\begin{equation}
\begin{split}
 & Pr\{\bX_{h_j, m_{h_j}} \leq \theta_{h_j} - C_{t, m_{h_j}} \} \nonumber \\
 & = Pr\{m_{h_j} \bX_{h_j, m_{h_j}} \leq m_{h_j} \theta_{h_j} - m_{h_j} C_{t, m_{h_j}} \}\\
 & \leq e^{-2\cdot \frac{1}{m_{h_ij}} \cdot (m_{h_j})^2 \cdot \frac{ (L+1) \ln t }{ m_{h_j} } } \\
 & = e^{-2(L+1) \ln t}\\
 & = t^{-2(L+1)}.\\
\end{split}
\end{equation}
Thus,
\begin{equation}
\begin{split}
 & Pr\{ \sum\limits_{j = 1}^{|\xA_{\xa*}|} a_{h_j}^* \bX_{h_j, m_{h_j}} \leq  R^*
   - \sum\limits_{j = 1}^{|\xA_{\xa*}|} a_{h_j}^* C_{t, m_{h_j}} \} \\
   & \leq |\xA_{\xa*}|
 t^{-2(L+1)}\\
 & \leq \quad L t^{-2(L+1)}.
\end{split}
\end{equation}
Similarly, we can get the upper bound of the probability for
inequality (\ref{equ:p2}):
\begin{equation}
\begin{split}
 & Pr\{ \sum\limits_{j = 1}^{|\xA_{\xa(t)}|} a_{p_j}(t) \bX_{p_j, m_{p_j}}
 \geq R_{\xa(t)} + \sum\limits_{j = 1}^{|\xA_{\xa(t)}|} a_{p_j}(t) C_{t, m_{p_j}}
 \} \\
 & \qquad \leq L t^{-2(L+1)}.
\end{split}
\end{equation}

Note that for $l \geq \left\lceil \frac{4(L+1) \ln n}{
\left(\frac{\Delta_{\xa(t)}}{L a_{\max}} \right)^2 } \right\rceil$,
\begin{equation}
 \label{equ:p4}
\begin{split}
 & R^* - R_{\xa(t)}  - 2 \sum\limits_{j = 1}^{|\xA_{\xa(t)}|} a_{p_j}(t) C_{t, m_{p_j}}  \\
 & = R^* - R_{\xa(t)} - 2 \sum\limits_{j = 1}^{|\xA_{\xa(t)}|} a_{p_j} \sqrt{ \frac{ (L+1) \ln t }{ m_{p_j}} } \\
 & \geq R^* - R_{\xa(t)} - L a_{\max} \sqrt{ \frac{ 4(L+1) \ln n }{ l } } \\
 & \geq R^* - R_{\xa(t)} - L a_{\max} \sqrt{ \frac{ 4(L+1) \ln n }{ 4(L+1) \ln n } \left( \frac{\Delta_{\xa(t)} }{L a_{\max}} \right)^2 } \\
 & \geq R^* - R_{\xa(t)} - \Delta_{\xa(t)} = 0.
\end{split}
\end{equation}

Equation (\ref{equ:p4}) implies that condition (\ref{equ:p1}) is
false when $l = \left\lceil \frac{4(L+1) \ln n}{
\left(\frac{\Delta_{\xa(t)}}{L a_{\max}} \right)^2 } \right\rceil$.
If we let $l = \left\lceil \frac{4(L+1) \ln n}{
\left(\frac{\Delta_{\min}}{L a_{\max}} \right)^2 } \right\rceil$,
then (\ref{equ:p1}) is false for all $\xa(t)$.

Therefore,
\begin{equation}
 \label{equ:p5}
\begin{split}
 & \mathbb{E}[\xT(n)]  \leq \left\lceil \frac{4(L+1) \ln n}{
\left(\frac{\Delta_{\min}}{L a_{\max}} \right)^2 } \right\rceil \\
 & + \sum\limits_{t = 1}^{\infty} \left( \sum\limits_{m_{h_1} = 1}^{t} \dots \sum\limits_{m_{h_{|\xA^*|}} = 1}^{t}
 \sum\limits_{m_{p_1} = l}^{t} \dots \sum\limits_{m_{p_{|\xA_{\xa(t)}|}} = l}^{t}  2 L t^{-2(L+1)} \right) \\
 & \leq  \frac{ 4 a_{\max}^2 L^2 (L+1)  \ln n }{ \left( \Delta_{\min} \right)^2 } + 1 + L \sum\limits_{t = 1}^\infty 2 t^{-2} \\
 & \leq  \frac{ 4 a_{\max}^2 L^2 (L+1)  \ln n }{ \left( \Delta_{\min} \right)^2 } + 1 + \frac{\pi^2}{3} L. \\
\end{split}
\end{equation}
So under LLR policy, we have:
\begin{equation}
 \label{equ:p6}
\begin{split}
 \mathfrak{R}^\pi_n (\Theta) & = R^* n - \mathbb{E}^\pi[ \sum \limits_{t = 1}^n R_{\pi(t)}(t) ]\\
 & = \sum\limits_{\xa: R_{\xa} < R^*} \Delta_{\xa} \mathbb{E}[T_{\xa}(n)] \\
 & \leq \Delta_{\max} \sum\limits_{\xa: R_{\xa} < R^*} \mathbb{E}[T_{\xa}(n)] \\
 & = \Delta_{\max} \sum\limits_{i = 1}^N  \mathbb{E}[\xT(n)] \\
 & \leq \left[\sum\limits_{i = 1}^N  \frac{ 4 a_{\max}^2 L^2 (L+1)  \ln n }{ \left( \Delta_{\min} \right)^2 } + N + \frac{\pi^2}{3} LN \right] \Delta_{\max}
 \\
 & \leq \left[\frac{ 4 a_{\max}^2 L^2 (L+1) N  \ln n }{ \left( \Delta_{\min} \right)^2 } + N + \frac{\pi^2}{3} LN \right] \Delta_{\max}. \\
\end{split}
\end{equation}
\end{IEEEproof}

\begin{remark} Note that when the set of action vectors consists of binary vectors
with a single ``1'', the problem formulation reduces to an
multi-armed bandit problem with $N$ independent arms. In this
special case, the LLR algorithm is equivalent to UCB1
in~\cite{Auer:2002}. Thus, our results generalize that prior work.
\end{remark}

\begin{remark} We have presented $\xF$ as a finite set in our problem
formation. We note that the LLR policy we have described and its
analysis actually also work with a more general formulation when
$\xF$ is an infinite set with the following additional constraints:
the maximization problem in (\ref{equ:choosemax}) always has at
least one solution; $\Delta_{\min}$ exists; $a_i$ is bounded. With
the above constraints, Algorithm \ref{alg:general} will work the
same and the conclusion and all the details of the proof of Theorem
\ref{c1:upperbound} can remain the same. \end{remark}

\begin{remark} Theorem \ref{c1:upperbound}
also holds for random variables $X_i, 1 \leq i \leq N$ that are not
i.i.d. over time, but with the only weaker assumption that
$E[X_i(t)|X_i(1), \ldots, X_i(t-1)] = \theta_i, \forall 1 \leq i
\leq N$. This is because the Chernoff-Hoeffding bound only needs a
weak assumption $E[X_i(t)|X_i(1), \ldots, X_i(t-1)] = \theta_i,
\forall 1 \leq i \leq N$.
\end{remark}

\section{Applications}\label{sec:example}

We now describe some applications and extensions of the LLR policy
for combinatorial network optimization in graphs where the edge
weights are unknown random variables.

\subsection{Maximum Weighted Matching}

Maximum Weighted Matching (MWM) problems are widely used in the many
optimization problems in wireless networks such as the prior work in
\cite{Balakrishnan:2004, Brzezinski:2006}. Given any graph $G = (V,
E)$, there is a weight associated with each edge and the objective
is to maximize the sum weights of a matching among all the matchings
in a given constraint set, i.e., the general formulation for MWM
problem is
\begin{equation}
 \begin{split}
  \max & \quad R_{\xa}^{MWM}  =  \sum\limits_{i = 1}^{|E|} a_{i} W_{i} \\
  s.t. & \quad \xa \text{ is a matching}
 \end{split}
\end{equation}
where $W_i$ is the weight associated with each edge $i$.

In many practical applications, the weights are unknown random
variables and we need to learn by selecting different matchings over
time. This kind of problem fits the general framework of our
proposed policy regarding the reward as the sum weight and a
matching as an arm. Our proposed LLR policy is a solution with
linear storage, and the regret polynomial in the number of edges,
and logarithmic in time.

Since there are various algorithms to solve the different variations
in the maximum weighted matching problems, such as the Hungarian
algorithm for the maximum weighted bipartite
matching~\cite{Kuhn:1955}, Edmonds's matching
algorithm~\cite{Edmonds:1965} for a general maximum matching. In
these cases, the computation time is also polynomial.

Here we present a general problem of multiuser channel allocations
in cognitive radio network. There are $M$ secondary users and $Q$
orthogonal channels. Each secondary user requires a single channel
for operation that does not conflict with the channels assigned to
the other users. Due to geographic dispersion, each secondary user
can potentially see different primary user occupancy behavior on
each channel. Time is divided into discrete decision rounds. The
throughput obtainable from spectrum opportunities on each
user-channel combination over a decision period is denoted as $S_{i,
j}$ and modeled as an arbitrarily-distributed random variable with
bounded support but unknown mean, i.i.d. over time. This random
process is assumed to have a mean $\theta_{i,j}$ that is unknown to
the users. The objective is to search for an allocation of channels
for all users that maximizes the expected sum throughput.

Assuming an interference model whereby at most one secondary user
can derive benefit from any channel, if the number of channels is
greater than the number of users, an optimal channel allocation
employs a one-to-one matching of users to channels, such that the
expected sum-throughput is maximized.

\begin{figure}[h]
  \centering
  \includegraphics[width=0.35\textwidth]{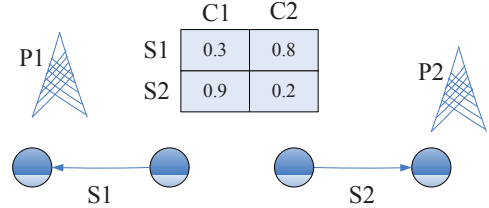}
 \caption{An illustrative scenario} \label{fig:1}
\end{figure}

Figure~\ref{fig:1} illustrates a simple scenario. There are two
secondary users (i.e., links) S1 and S2, that are each assumed to be
in interference range of each other. S1 is proximate to primary user
P1 who is operating on channel 1. S2 is proximate to primary user P2
who is operating on channel 2. The matrix shows the corresponding
$\Theta$, i.e., the throughput each secondary user could derive from
being on the corresponding channel. In this simple example, the
optimal matching is for secondary user 1 to be allocated channel 2
and user 2 to be allocated channel 1. Note, however, that, in our
formulation, the users are not \emph{a priori} aware of the matrix
of mean values, and therefore must follow a sequential learning
policy.

Note that this problem can be formulated as a multi-armed bandits
with linear regret, in which each arm corresponds to a matching of
the users to channels, and the reward corresponds to the
sum-throughput. In this channel allocation problem, there is $M
\times Q$ unknown random variables, and the number of arms are
$P(Q,M)$, which can grow exponentially in the number of unknown
random variables. Following the convention, instead of denoting the
variables as a vector, we refer it as a $M$ by $Q$ matrix. So the
reward as each time slot by choosing a permutation $\xa$ is
expressed as:
\begin{equation}
 \begin{split}
  R_{\xa}  = & \sum\limits_{i = 1}^M \sum\limits_{j = 1}^Q a_{i,j}
  S_{i,j} \\
 \end{split}
\end{equation}
where $\xa \in \xF$, $\xF$ is a set with all permutations, which is
defined as:
\begin{equation}
 \begin{split}
  \xF & = \{\xa: a_{i,j} \in \{0, 1\}, \forall i,j \wedge \sum\limits_{i = 1}^Q a_{i,j} = 1 \wedge \sum\limits_{j = 1}^Q a_{i,j} = 1 \}. \\
  %& \sum\limits_{i = 1}^N a_{i,j} = 1 \\
  %& \sum\limits_{j = 1}^N a_{i,j} = 1
 \end{split}
\end{equation}

We use two $M$ by $Q$ matrices to store the information after we
play an arm at each time slot. One is $(\hat{ \theta }_{i,j})_{M
\times Q}$ in which $\hat{ \theta }_{i,j}$ is the average (sample
mean) of all the observed values of channel $j$ by user $i$ up to
the current time slot (obtained through potentially different sets
of arms over time). The other one is $(m_{i,j})_{M \times Q}$ in
which $m_{i,j}$ is the number of times that channel $j$ has been
observed by user $i$ up to the current time slot.

Applying Algorithm \ref{alg:general},  we get a linear storage
policy for which $(\hat{ \theta }_{i,j})_{M \times Q}$ and
$(m_{i,j})_{M \times Q}$ are stored and updated at each time slot.
The regret is polynomial in the number of users and channels, and
logarithmic in time. Also, the computation time for the policy is
also polynomial since (\ref{equ:choosemax}) in Algorithm
\ref{alg:general} now becomes the following deterministic maximum
weighted bipartite matching problem
\begin{equation}
    \label{equ:channelallocation}
    \arg \max\limits_{\xa \in \xF} \sum\limits_{(i,j) \in \xA_\xa} \left( \hat{ \theta }_{i,j} + \sqrt{ \frac{ (L+1) \ln n }{ m_{i,j}}
    } \right)
\end{equation}
on the bipartite graph of users and channels with edge weights
$\left(\hat{ \theta }_{i,j} + \sqrt{ \frac{ (L+1) \ln n }{ m_{i,j}}
}\right)$. It could be solved with polynomial computation time
(e.g., using the Hungarian algorithm~\cite{Kuhn:1955}). Note that $L
= \max\limits_\xa |\xA_\xa| = \min\{M, Q\}$ for this problem, which
is less than $M\times Q$ so that the bound of regret is tighter. The
regret is $O(\min\{M, Q\}^3 M Q\log n)$ following Theorem
\ref{c1:upperbound}.

\subsection{Shortest Path}

Shortest Path (SP) problem is another example where the underlying
deterministic optimization can be done with polynomial computation
time. If the given directed graph is denoted as $G=(V, E)$ with the
source node $s$ and the destination node $d$, and the cost (e.g.,
the transmission delay) associated with edge $(i,j)$ is denoted as
$D_{i,j} \geq 0$, the objective is find the path from $s$ to $d$
with the minimum sum cost, i.e.,
\begin{align}
  \min & \quad C_{\xa}^{SP} = \sum\limits_{(i,j) \in E} a_{i,j} D_{i,j} \\
  s.t. & \quad a_{i,j} \in \{0, 1 \}, \forall (i, j) \in E \label{equ:23}\\
  & \quad \forall i, \sum\limits_j a_{i,j} - \sum\limits_j a_{j,i} =
  \left\{ \begin{array}
 {l@{\quad:\quad}l}
 1 &  i = s\\\label{equ:24}
 -1 & i = t\\
 0 & \text{otherwise}
\end{array} \right.
\end{align}
where equation (\ref{equ:23}) and (\ref{equ:24}) defines a feasible
set $\xF$, such that $\xF$ is the set of all possible pathes from
$s$ to $d$. When $(D_{ij})$ are random variables with bounded
support but unknown mean, i.i.d. over time, an dynamic learning
policy is needed for this multi-armed bandit formulation.

Note that corresponding to the LLR policy with the objective to
maximize the rewards, a direct variation of it is to find the
minimum linear cost defined on finite constraint set $\xF$, by
changing the maximization problem in to a minimization problem. For
clarity, this straightforward modification of LLR is shown below in
Algorithm \ref{alg:cost}, which we refer to as Learning with Linear
Costs (LLC).

%%%%%%%%%%%%%%%%%%%%%%
%%%%%%%%%%%%%%%%%%%%%%%
\begin{algorithm} [ht]
\caption{Learning with Linear Cost (LLC)} \label{alg:cost}

\begin{algorithmic}[1]
%----------------------------------------------------------------
\State \commnt{Initialization part is same as in Algorithm
\ref{alg:general}}

\State \commnt{Main loop}

\While {1}
    \State $n = n + 1$;
    \State Play an arm $\xa$ which solves the minimization problem
    \begin{equation}
    \label{equ:choosemin}
    \xa = \arg\min\limits_{\xa \in \xF} \sum\limits_{i \in \xA_\xa} a_i \left( \hat{ \theta }_{i} - \sqrt{ \frac{ (L+1) \ln n }{ m_{i}}
    } \right);
    \end{equation}
    \State Update $(\hat{\theta}_{i})_{1 \times N}$, $(m_{i})_{1 \times N}$
    accordingly;
\EndWhile
\end{algorithmic}
\end{algorithm}
%%%%%%%%%%%%%%%%%%%%%%%
%%%%%%%%%%%%%%%%%%%%%%%

LLC (Algorithm \ref{alg:cost}) is a policy for a general multi-armed
bandit problem with linear cost defined on any constraint set. It is
directly derived from the LLR policy (Algorithm \ref{alg:general}),
so Theorem \ref{c1:upperbound} also holds for LLC, where the regret
is defined as:
\begin{equation}
 \mathfrak{R}^\pi_n (\Theta) = E^\pi[ \sum \limits_{t = 1}^n C_{\pi(t)}(t) ] -n
 C^*
\end{equation}
where $C^*$ represents the minimum cost, which is cost of the
optimal arm.

Using the LLC policy, we map each path between $s$ and $t$ as an
arm. The number of unknown variables are $|E|$, while the number of
arms could grow exponentially in the worst case. Since there exist
polynomial computation time algorithms such as Dijkstra's
algorithm~\cite{Dijkstra:1959} and Bellman-Ford
algorithm~\cite{Bellman:1958,Ford:1956} for the shortest path
problem, we could apply these algorithms to solve
(\ref{equ:choosemin}) with edge cost $\hat{ \theta }_{i} - \sqrt{
\frac{ (L+1) \ln n }{ m_{i}}}$. LLC is thus an efficient policy to
solve the multi-armed bandit formulation of the shortest path
problem with linear storage, polynomial computation time. Note that
$L = \max\limits_\xa |\xA_\xa| = |E|$. Regret is $O(|E|^4 \log n)$.

Another related problem is the Shortest Path Tree (SPT), where
problem formulation is similar, and the objective is to find a
subgraph of the given graph with the minimum total cost between a
selected root $s$ node and all other nodes. It is expressed
as~\cite{Karup:2004,Bazaraa:2009}:
\begin{align}
  \min & \quad C_{\xa}^{SPT} = \sum\limits_{(i,j) \in E} a_{i,j} D_{i,j} \\
  s.t. & \quad a_{i,j} \in \{0, 1 \}, \forall (i, j) \in E \label{equ:28}\\
  & \quad \sum\limits_{(j,i) \in \mathcal{BS}(i)} a_{j,i} - \sum\limits_{(i,j) \in \mathcal{FS}(i)}
  a_{i,j} \nonumber
  \\ & \quad\quad
  =
  \left\{ \begin{array}
 {l@{\quad:\quad}l}
 -n+1 &  i = s \\
 1 & i \in V/\{s\} \label{equ:27}\\
\end{array}
\right.
\end{align}
where $\mathcal{BS}(i) = \{(u,v) \in E: v=i \}$, $\mathcal{FS}(i) =
\{(u,v) \in E: u=i \}$. (\ref{equ:27}) and (\ref{equ:28}) defines
the constraint set $\xF$. We can also use the polynomial computation
time algorithms such as Dijkstra's algorithm and Bellman-Ford
algorithm to solve (\ref{equ:choosemin}) for the LLC policy.

\subsection{Minimum Spanning Tree}

Minimum Spanning Tree (MST) is another combinatorial optimization
with polynomial computation time algorithms, such as Prim's
algorithm \cite{Prim:1957} and Kruskal's algorithm
\cite{Kruskal:1956}. The objective for the MST problem can be simply
presented as
\begin{equation}
 \min\limits_{\xa \in \xF} C_{\xa}^{MST} = \sum\limits_{(i,j) \in E} a_{i,j} D_{i,j}
\end{equation}
where $\xF$ is the set of all spanning trees in the graph.

With the LLC policy, each spanning tree is treated as an arm, and $L
= |E|$. Regret bound also grows as $O(|E|^4 \log n)$.

\section{Numerical Simulation Results} \label{sec:simulation}

\begin{figure}[h]
  \centering
  \includegraphics[width=0.48\textwidth]{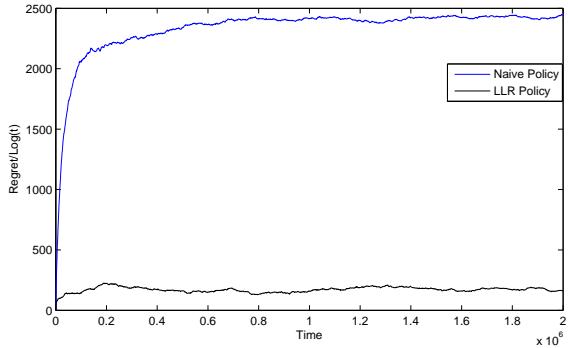}
 \caption{Simulation Results of a system with $7$ orthogonal channels and $4$ users.} \label{fig:7by4}
\end{figure}

\begin{figure}[h]
  \centering
  \includegraphics[width=0.48\textwidth]{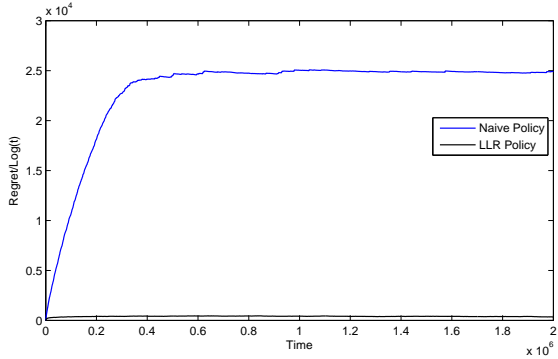}
 \caption{Simulation Results of a system with $9$ orthogonal channels and $5$ users.} \label{fig:9by5}
\end{figure}

We present in the section the numerical simulation results with the
example of multiuser channel allocations in cognitive radio network.

Fig \ref{fig:7by4} shows the simulation results of using LLR policy
compared with the naive policy in \ref{sec:naive}. We assume that
the system consists of $Q = 7$ orthogonal channels in and $M = 4$
secondary users. The throughput $\{S_{i,j}(t)\}_{t \geq 1}$ for the
user-channel combination is an i.i.d. Bernoulli  process with mean
$\theta_{i,j}$ ($(\theta_{i,j})$ is unknown to the players) shown as
below:
\begin{equation}
(\theta_{i,j})=\left(
\begin{array}{ccccccc}
 0.3 & 0.5 & \fbox{\parbox{0.4cm}{0.9}} & 0.7 & 0.8 & 0.9 & 0.6\\
 0.2 & 0.2 & 0.3 & 0.4 & \fbox{\parbox{0.4cm}{0.5}} & 0.4 & 0.5\\
 \fbox{\parbox{0.4cm}{0.8}} & 0.6 & 0.5 & 0.4 & 0.7 & 0.2 & 0.8\\
 0.9 & 0.2 & 0.2 & 0.8 & 0.3 & \fbox{\parbox{0.4cm}{0.9}} & 0.6
\end{array}
\right)
\end{equation}
where the components in the box are in the optimal arm. Note that
$P(7,4) = 840$ while $7\times4 = 28$, so the storage used for the
naive approach is $30$ times more than the LLR policy. Fig
\ref{fig:7by4} shows the regret (normalized with respect to the
logarithm of time) over time for the naive policy and the LLR
policy. We can see that under both policies the regret grows
logarithmically in time. But the regret for the naive policy is a
lot higher than that of the LLR policy.

Fig \ref{fig:9by5} is another example of the case when $Q = 9$ and
$M = 5$. The throughput is also assumed to be an i.i.d. Bernoulli
process, with the following mean:
\begin{equation}
\begin{split}
  & (\theta_{i,j}) = \\
  & \left(
\begin{array}{ccccccccc}
 0.3 & 0.5 & \fbox{\parbox{0.4cm}{0.9}} & 0.7 & 0.8 & 0.9 & 0.6 & 0.8 & 0.7\\
 0.2 & 0.2 & 0.3 & 0.4 & 0.5 & 0.4 & 0.5 & 0.6 & \fbox{\parbox{0.4cm}{0.9}}\\
 0.8 & 0.6 & 0.5 & 0.4 & 0.7 & 0.2 & \fbox{\parbox{0.4cm}{0.8}} & 0.2 & 0.8\\
 \fbox{\parbox{0.4cm}{0.9}} & 0.2 & 0.2 & 0.8 & 0.3 & 0.9 & 0.6 & 0.5 & 0.4\\
 0.6 & 0.7 & 0.5 & 0.7 & 0.6 & \fbox{\parbox{0.4cm}{0.8}} & 0.2 & 0.6 & 0.8
\end{array}
\right).
\end{split}
\end{equation}

For this example, $P(9,5) = 15120$, which is much higher than
$9\times5 = 45$ (about $336$ times higher), so the storage used by
the naive policy grows much faster than the LLR policy. Comparing
with the regrets shown in Table \ref{tab:regret} for both examples
when $t = 2\times10^6$, we can see that the regret also grows much
faster for the naive policy.

\begin{table}[!h]
\tabcolsep 0pt \normalsize \caption{Regret when $t = 2\times10^6$}
\label{tab:regret} \vspace*{-12pt}
\begin{center}
\def\temptablewidth{0.4\textwidth}
{\rule{\temptablewidth}{1pt}}
\begin{tabular*}{0.4\textwidth}{@{\extracolsep{\fill}}c|cc}
 & Naive Policy & LLR  \\   \hline
 7 channels, 4 users \; & 2443.6 & 163.6\\ \hline
 9 channels, 5 users \; & 24892.6 & 345.2
 \end{tabular*}
{\rule{\temptablewidth}{1pt}}
\end{center}
\end{table}

\section{K Simultaneous Actions} \label{sec:klargest}

The reward-maximizing LLR policy presented in Algorithm
\ref{alg:general} and the corresponding cost-minimizing LLC policy
presented in \ref{alg:cost} can also be extended to the setting
where $K$ arms are played at each time slot. The goal is to maximize
the total rewards (or minimize the total costs) obtained by these
$K$ arms. For brevity, we only present the policy for the
reward-maximization problem; the extension to cost-minimization is
straightforward. The modified LLR-K policy for picking the $K$ best
arms are shown in Algorithm \ref{alg:karm}.
%%%%%%%%%%%%%%%%%%%%%%
%%%%%%%%%%%%%%%%%%%%%%%
\begin{algorithm} [ht]
\caption{Learning with Linear Rewards while selecting $K$ arms
(LLR-K)} \label{alg:karm}

\begin{algorithmic}[1]
%----------------------------------------------------------------
\State \commnt{Initialization part is same as in Algorithm
\ref{alg:general}}

\State \commnt{Main loop}

\While {1}
    \State $n = n + 1$;
    \State Play arms $\{\xa\}_{K} \in \xF$ with $K$ largest values in (\ref{equ:k})
    \begin{equation}
    \label{equ:k}
      \sum\limits_{i \in \xA_\xa} a_i \left( \hat{ \theta }_{i} + \sqrt{ \frac{ (L+1) \ln n }{ m_{i}}
    } \right);
    \end{equation}
    \State Update $(\hat{\theta}_{i})_{1 \times N}$, $(m_{i})_{1 \times N}$ for all arms
    accordingly;
\EndWhile
\end{algorithmic}
\end{algorithm}
%%%%%%%%%%%%%%%%%%%%%%%
%%%%%%%%%%%%%%%%%%%%%%%

Theorem \ref{c1:maxk} states the upper bound of the regret for the
extended LLR-K policy.

\theorem\label{c1:maxk} The expected regret under the LLR-K policy
with $K$ arms selection is at most
\begin{equation}
\left[\frac{ 4 a_{\max}^2 L^2 (L+1) N \ln n }{ \left( \Delta_{\min}
\right)^2 } + N + \frac{\pi^2}{3} L K^{2L} N \right] \Delta_{\max}.
\end{equation}

\begin{IEEEproof}

The proof is similar to the proof of Theorem \ref{c1:upperbound},
but now we have a set of $K$ arms with $K$ largest expected rewards
as the optimal arms. We denote this set as $\mathfrak{A^*} = \{
\xa^{*, k}, 1 \leq k \leq K \}$ where $\xa^{*, k}$ is the arm with
$k$-th largest expected reward. As in the proof of Theorem
\ref{c1:upperbound}, we define $\xT(n)$ as a counter when a
non-optimal arm is played in the same way. Equation (\ref{equ:f1}),
(\ref{equ:f2}), (\ref{equ:10}) and (\ref{equ:13}) still hold.

Note that each time when $\xI(t) = 1$, there exists some arm such
that a non-optimal arm is picked for which $m_{i}$ is the minimum in
this arm. We denote this arm as $\xa(t)$. Note that $\xa(t)$ means
there exists $m$, $1\leq m \leq K$, such that the following holds:
\begin{equation}
\begin{split}
 \xT(n)
 & \leq l + \sum\limits_{t = N}^n \{ \sum\limits_{j \in \xA_{\xa^{*,m}}} a_j^{*,m} (\bX_{j, m_j(t)} + C_{t, m_j(t)} ) \\
 & \leq \sum\limits_{j \in \xA_{\xa(t)}} a_j(t) (\bX_{j, m_j(t)} + C_{t, m_j(t)} ), \xT(t) \geq l
 \}.
\end{split}
\end{equation}

Since at each time $K$ arms are played, so at time $t$, an random
variable could be observed up to $Kt$ times. Then (\ref{equ:12})
should be modified as:

\begin{equation}
\begin{split}
 \xT(n)
 & \leq l + \sum\limits_{t = 1}^{\infty} \sum\limits_{m_{h_1} = 1}^{Kt} \dots \sum\limits_{m_{h_{|\xA^{*,m}|}} = 1}^{Kt}
 \sum\limits_{m_{p_1} = l}^{Kt} \dots \sum\limits_{m_{p_{|\xA_{\xa(t)}|}} = l}^{Kt}    \\
 & \quad\quad \{\sum\limits_{j = 1}^{|\xA_{\xa^{*,m}}|} a_{h_j}^{*,m} (\bX_{h_j, m_{h_j}} + C_{t, m_{h_j}} ) \\
 & \quad\quad \leq \sum\limits_{j = 1}^{|\xA_{\xa(t)}|} a_{p_j}(t) (\bX_{p_j, m_{p_j}} + C_{t, m_{p_j}} )
 \}.
\end{split}
\end{equation}

Equation (\ref{equ:p1}) to (\ref{equ:p4}) are similar by
substituting $\xa^*$ with $\xa^{*,m}$. So, we have:
\begin{equation}
 \label{equ:p5}
\begin{split}
 & \mathbb{E}[\xT(n)]  \leq \left\lceil \frac{4(L+1) \ln n}{
\left(\frac{\Delta_{\min}}{L a_{\max}} \right)^2 } \right\rceil \\
 & + \sum\limits_{t = 1}^{\infty} \left( \sum\limits_{m_{h_1} = 1}^{Kt} \dots \sum\limits_{m_{h_{|\xA^*|}} = 1}^{Kt}
 \sum\limits_{m_{p_1} = l}^{Kt} \dots \sum\limits_{m_{p_{|\xA_{\xa(t)}|}} = l}^{Kt}  2 L t^{-2(L+1)} \right) \\
 & \leq  \frac{ 4 a_{\max}^2 L^2 (L+1)  \ln n }{ \left( \Delta_{\min} \right)^2 } + 1 + \frac{\pi^2}{3} L K^{2L}. \\
\end{split}
\end{equation}
Hence, we get the upper bound for the regret as:
\begin{equation}
 \label{equ:p6}
\begin{split}
 \mathfrak{R}^\pi_n (\Theta)
 & \leq \left[\frac{ 4 a_{\max}^2 L^2 (L+1) N  \ln n }{ \left( \Delta_{\min} \right)^2 } + N + \frac{\pi^2}{3} L K^{2L} N \right] \Delta_{\max}. \\
\end{split}
\end{equation}

\end{IEEEproof}

\section{Conclusion}\label{sec:conclusion}

We have considered multi-armed bandit problems that provide for arms
with rewards that are a linear function of a smaller set of random
variables with unknown means. For such problems, if the number of
arms is exponentially large in the number of underlying random
variables, existing arm-based index policies such as the well-known
UCB1~\cite{Auer:2002} have poor performance in terms of storage,
computation, and regret. The LLR and LLR policies we have presented
are smarter in that they store and make decisions at each time based
on the stochastic observations of the underlying unknown-mean random
variables alone; they require only linear storage and result in a
regret that is bounded by a polynomial function of the number of
unknown-mean random variables. If the deterministic version of the
corresponding combinatorial optimization problem can be solved in
polynomial time, our policy will also require only polynomial
computation per step. We have shown a number of problems in the
context of networks where this formulation would be useful,
including maximum-weight matching, shortest path and spanning tree
computations.

While this work has provided useful insights into real-world linear
combinatorial optimization with unknown-mean random coefficients,
there are many interesting open problems to be explored in the
future. One open question is to derive a lower bound on the regret
achievable by any policy for this problem. We conjecture on
intuitive grounds that it is not possible to have regret lower than
$\Omega(N \log n)$, but this remains to be proved rigorously. It is
unclear whether the lower bound can be any higher than this, and
hence, it is unclear whether it is possible to prove an upper bound
on regret for some policy that is better than the $O(N^4 \log n)$
upper bound shown in our work.

In the context of channel access in cognitive radio networks, other
researchers have recently developed distributed policies in which
different users each select an arm independently~\cite{Liu:Zhao,
Anandkumar:2010}. A closely related problem in this setting would be
to have distributed users selecting different elements of the action
vector independently. The design and analysis of such distributed
policies is an open problem.

Finally, it would be of great interest to see if it is possible to
also tackle non-linear reward functions, at least in structured
cases that have proved to be tractable in deterministic settings,
such as convex functions.% or sub-modular set functions. (Daniel

\end{document}